\nonstopmode
\input amstex
\input amsppt.sty   
%
\catcode`\@=11

\def\East#1#2{\setboxz@h{$\m@th\ssize\;{#1}\;\;$}%
 \setbox@ne\hbox{$\m@th\ssize\;{#2}\;\;$}\setbox\tw@\hbox{$\m@th#2$}%
 \dimen@\minaw@
 \ifdim\wdz@>\dimen@ \dimen@\wdz@ \fi  \ifdim\wd@ne>\dimen@ \dimen@\wd@ne \fi
 \ifdim\wd\tw@>\z@
  \mathrel{\mathop{\hbox to\dimen@{\rightarrowfill}}\limits^{#1}_{#2}}%
 \else
  \mathrel{\mathop{\hbox to\dimen@{\rightarrowfill}}\limits^{#1}}%
 \fi}
\def\West#1#2{\setboxz@h{$\m@th\ssize\;\;{#1}\;$}%
 \setbox@ne\hbox{$\m@th\ssize\;\;{#2}\;$}\setbox\tw@\hbox{$\m@th#2$}%
 \dimen@\minaw@
 \ifdim\wdz@>\dimen@ \dimen@\wdz@ \fi \ifdim\wd@ne>\dimen@ \dimen@\wd@ne \fi
 \ifdim\wd\tw@>\z@
  \mathrel{\mathop{\hbox to\dimen@{\leftarrowfill}}\limits^{#1}_{#2}}%
 \else
  \mathrel{\mathop{\hbox to\dimen@{\leftarrowfill}}\limits^{#1}}%
 \fi}
\font\arrow@i=lams1
\font\arrow@ii=lams2
\font\arrow@iii=lams3
\font\arrow@iv=lams4
\font\arrow@v=lams5
\newbox\zer@
\newdimen\standardcgap
\standardcgap=40\p@
\newdimen\hunit
\hunit=\tw@\p@
\newdimen\standardrgap
\standardrgap=32\p@
\newdimen\vunit
\vunit=1.6\p@
\def\Cgaps#1{\RIfM@
  \standardcgap=#1\standardcgap\relax \hunit=#1\hunit\relax
 \else \nonmatherr@\Cgaps \fi}
\def\Rgaps#1{\RIfM@
  \standardrgap=#1\standardrgap\relax \vunit=#1\vunit\relax
 \else \nonmatherr@\Rgaps \fi}
\newdimen\getdim@
\def\getcgap@#1{\ifcase#1\or\getdim@\z@\else\getdim@\standardcgap\fi}
\def\getrgap@#1{\ifcase#1\getdim@\z@\else\getdim@\standardrgap\fi}
\def\cgaps#1{\RIfM@
 \cgaps@{#1}\edef\getcgap@##1{\i@=##1\relax\the\toks@}\toks@{}\else
 \nonmatherr@\cgaps\fi}
\def\rgaps#1{\RIfM@
 \rgaps@{#1}\edef\getrgap@##1{\i@=##1\relax\the\toks@}\toks@{}\else
 \nonmatherr@\rgaps\fi}
\def\Gaps@@{\gaps@@}
\def\cgaps@#1{\toks@{\ifcase\i@\or\getdim@=\z@}%
 \gaps@@\standardcgap#1;\gaps@@\gaps@@
 \edef\next@{\the\toks@\noexpand\else\noexpand\getdim@\noexpand\standardcgap
  \noexpand\fi}%
 \toks@=\expandafter{\next@}}
\def\rgaps@#1{\toks@{\ifcase\i@\getdim@=\z@}%
 \gaps@@\standardrgap#1;\gaps@@\gaps@@
 \edef\next@{\the\toks@\noexpand\else\noexpand\getdim@\noexpand\standardrgap
  \noexpand\fi}%
 \toks@=\expandafter{\next@}}
\def\gaps@@#1#2;#3{\mgaps@#1#2\mgaps@
 \edef\next@{\the\toks@\noexpand\or\noexpand\getdim@
  \noexpand#1\the\mgapstoks@@}%
 \global\toks@=\expandafter{\next@}%
 \DN@{#3}%
 \ifx\next@\Gaps@@\gdef\next@##1\gaps@@{}\else
  \gdef\next@{\gaps@@#1#3}\fi\next@}
\def\mgaps@#1{\let\mgapsnext@#1\FN@\mgaps@@}
\def\mgaps@@{\ifx\next\space@\DN@. {\FN@\mgaps@@}\else
 \DN@.{\FN@\mgaps@@@}\fi\next@.}
\def\mgaps@@@{\ifx\next\w\let\next@\mgaps@@@@\else
 \let\next@\mgaps@@@@@\fi\next@}
\newtoks\mgapstoks@@
\def\mgaps@@@@@#1\mgaps@{\getdim@\mgapsnext@\getdim@#1\getdim@
 \edef\next@{\noexpand\getdim@\the\getdim@}%
 \mgapstoks@@=\expandafter{\next@}}
\def\mgaps@@@@\w#1#2\mgaps@{\mgaps@@@@@#2\mgaps@
 \setbox\zer@\hbox{$\m@th\hskip15\p@\tsize@#1$}%
 \dimen@\wd\zer@
 \ifdim\dimen@>\getdim@ \getdim@\dimen@ \fi
 \edef\next@{\noexpand\getdim@\the\getdim@}%
 \mgapstoks@@=\expandafter{\next@}}
\def\changewidth#1#2{\setbox\zer@\hbox{$\m@th#2}%
 \hbox to\wd\zer@{\hss$\m@th#1$\hss}}
\atdef@({\FN@\ARROW@}
\def\ARROW@{\ifx\next)\let\next@\OPTIONS@\else
 \DN@{\csname\string @(\endcsname}\fi\next@}
\newif\ifoptions@
\def\OPTIONS@){\ifoptions@\let\next@\relax\else
 \DN@{\options@true\begingroup\optioncodes@}\fi\next@}
\newif\ifN@
\newif\ifE@
\newif\ifNESW@
\newif\ifH@
\newif\ifV@
\newif\ifHshort@
\expandafter\def\csname\string @(\endcsname #1,#2){%
 \ifoptions@\let\next@\endgroup\else\let\next@\relax\fi\next@
 \N@false\E@false\H@false\V@false\Hshort@false
 \ifnum#1>\z@\E@true\fi
 \ifnum#1=\z@\V@true\tX@false\tY@false\a@false\fi
 \ifnum#2>\z@\N@true\fi
 \ifnum#2=\z@\H@true\tX@false\tY@false\a@false\ifshort@\Hshort@true\fi\fi
 \NESW@false
 \ifN@\ifE@\NESW@true\fi\else\ifE@\else\NESW@true\fi\fi
 \arrow@{#1}{#2}%
 \global\options@false
 \global\scount@\z@\global\tcount@\z@\global\arrcount@\z@
 \global\s@false\global\sxdimen@\z@\global\sydimen@\z@
 \global\tX@false\global\tXdimen@i\z@\global\tXdimen@ii\z@
 \global\tY@false\global\tYdimen@i\z@\global\tYdimen@ii\z@
 \global\a@false\global\exacount@\z@
 \global\x@false\global\xdimen@\z@
 \global\X@false\global\Xdimen@\z@
 \global\y@false\global\ydimen@\z@
 \global\Y@false\global\Ydimen@\z@
 \global\p@false\global\pdimen@\z@
 \global\label@ifalse\global\label@iifalse
 \global\dl@ifalse\global\ldimen@i\z@
 \global\dl@iifalse\global\ldimen@ii\z@
 \global\short@false\global\unshort@false}
\newif\iflabel@i
\newif\iflabel@ii
\newcount\scount@
\newcount\tcount@
\newcount\arrcount@
\newif\ifs@
\newdimen\sxdimen@
\newdimen\sydimen@
\newif\iftX@
\newdimen\tXdimen@i
\newdimen\tXdimen@ii
\newif\iftY@
\newdimen\tYdimen@i
\newdimen\tYdimen@ii
\newif\ifa@
\newcount\exacount@
\newif\ifx@
\newdimen\xdimen@
\newif\ifX@
\newdimen\Xdimen@
\newif\ify@
\newdimen\ydimen@
\newif\ifY@
\newdimen\Ydimen@
\newif\ifp@
\newdimen\pdimen@
\newif\ifdl@i
\newif\ifdl@ii
\newdimen\ldimen@i
\newdimen\ldimen@ii
\newif\ifshort@
\newif\ifunshort@
\def\zero@#1{\ifnum\scount@=\z@
 \if#1e\global\scount@\m@ne\else
 \if#1t\global\scount@\tw@\else
 \if#1h\global\scount@\thr@@\else
 \if#1'\global\scount@6 \else
 \if#1`\global\scount@7 \else
 \if#1(\global\scount@8 \else
 \if#1)\global\scount@9 \else
 \if#1s\global\scount@12 \else
 \if#1H\global\scount@13 \else
 \Err@{\Invalid@@ option \string\0}\fi\fi\fi\fi\fi\fi\fi\fi\fi
 \fi}
\def\one@#1{\ifnum\tcount@=\z@
 \if#1e\global\tcount@\m@ne\else
 \if#1h\global\tcount@\tw@\else
 \if#1t\global\tcount@\thr@@\else
 \if#1'\global\tcount@4 \else
 \if#1`\global\tcount@5 \else
 \if#1(\global\tcount@10 \else
 \if#1)\global\tcount@11 \else
 \if#1s\global\tcount@12 \else
 \if#1H\global\tcount@13 \else
 \Err@{\Invalid@@ option \string\1}\fi\fi\fi\fi\fi\fi\fi\fi\fi
 \fi}
\def\a@#1{\ifnum\arrcount@=\z@
 \if#10\global\arrcount@\m@ne\else
 \if#1+\global\arrcount@\@ne\else
 \if#1-\global\arrcount@\tw@\else
 \if#1=\global\arrcount@\thr@@\else
 \Err@{\Invalid@@ option \string\a}\fi\fi\fi\fi
 \fi}
\def\ds@(#1;#2){\ifs@\else
 \global\s@true
 \sxdimen@\hunit \global\sxdimen@#1\sxdimen@\relax
 \sydimen@\vunit \global\sydimen@#2\sydimen@\relax
 \fi}
\def\dtX@(#1;#2){\iftX@\else
 \global\tX@true
 \tXdimen@i\hunit \global\tXdimen@i#1\tXdimen@i\relax
 \tXdimen@ii\vunit \global\tXdimen@ii#2\tXdimen@ii\relax
 \fi}
\def\dtY@(#1;#2){\iftY@\else
 \global\tY@true
 \tYdimen@i\hunit \global\tYdimen@i#1\tYdimen@i\relax
 \tYdimen@ii\vunit \global\tYdimen@ii#2\tYdimen@ii\relax
 \fi}
\def\da@#1{\ifa@\else\global\a@true\global\exacount@#1\relax\fi}
\def\dx@#1{\ifx@\else
 \global\x@true
 \xdimen@\hunit \global\xdimen@#1\xdimen@\relax
 \fi}
\def\dX@#1{\ifX@\else
 \global\X@true
 \Xdimen@\hunit \global\Xdimen@#1\Xdimen@\relax
 \fi}
\def\dy@#1{\ify@\else
 \global\y@true
 \ydimen@\vunit \global\ydimen@#1\ydimen@\relax
 \fi}
\def\dY@#1{\ifY@\else
 \global\Y@true
 \Ydimen@\vunit \global\Ydimen@#1\Ydimen@\relax
 \fi}
\def\p@@#1{\ifp@\else
 \global\p@true
 \pdimen@\hunit \divide\pdimen@\tw@ \global\pdimen@#1\pdimen@\relax
 \fi}
\def\L@#1{\iflabel@i\else
 \global\label@itrue \gdef\label@i{#1}%
 \fi}
\def\l@#1{\iflabel@ii\else
 \global\label@iitrue \gdef\label@ii{#1}%
 \fi}
\def\dL@#1{\ifdl@i\else
 \global\dl@itrue \ldimen@i\hunit \global\ldimen@i#1\ldimen@i\relax
 \fi}
\def\dl@#1{\ifdl@ii\else
 \global\dl@iitrue \ldimen@ii\hunit \global\ldimen@ii#1\ldimen@ii\relax
 \fi}
\def\s@{\ifunshort@\else\global\short@true\fi}
\def\uns@{\ifshort@\else\global\unshort@true\global\short@false\fi}
\def\optioncodes@{\let\0\zero@\let\1\one@\let\a\a@\let\ds\ds@\let\dtX\dtX@
 \let\dtY\dtY@\let\da\da@\let\dx\dx@\let\dX\dX@\let\dY\dY@\let\dy\dy@
 \let\p\p@@\let\L\L@\let\l\l@\let\dL\dL@\let\dl\dl@\let\s\s@\let\uns\uns@}
\def\slopes@{\\161\\152\\143\\134\\255\\126\\357\\238\\349\\45{10}\\56{11}%
 \\11{12}\\65{13}\\54{14}\\43{15}\\32{16}\\53{17}\\21{18}\\52{19}\\31{20}%
 \\41{21}\\51{22}\\61{23}}
\newcount\tan@i
\newcount\tan@ip
\newcount\tan@ii
\newcount\tan@iip
\newdimen\slope@i
\newdimen\slope@ip
\newdimen\slope@ii
\newdimen\slope@iip
\newcount\angcount@
\newcount\extracount@
\def\slope@{{\slope@i=\secondy@ \advance\slope@i-\firsty@
 \ifN@\else\multiply\slope@i\m@ne\fi
 \slope@ii=\secondx@ \advance\slope@ii-\firstx@
 \ifE@\else\multiply\slope@ii\m@ne\fi
 \ifdim\slope@ii<\z@
  \global\tan@i6 \global\tan@ii\@ne \global\angcount@23
 \else
  \dimen@\slope@i \multiply\dimen@6
  \ifdim\dimen@<\slope@ii
   \global\tan@i\@ne \global\tan@ii6 \global\angcount@\@ne
  \else
   \dimen@\slope@ii \multiply\dimen@6
   \ifdim\dimen@<\slope@i
    \global\tan@i6 \global\tan@ii\@ne \global\angcount@23
   \else
    \tan@ip\z@ \tan@iip \@ne
    \def\\##1##2##3{\global\angcount@=##3\relax
     \slope@ip\slope@i \slope@iip\slope@ii
     \multiply\slope@iip##1\relax \multiply\slope@ip##2\relax
     \ifdim\slope@iip<\slope@ip
      \global\tan@ip=##1\relax \global\tan@iip=##2\relax
     \else
      \global\tan@i=##1\relax \global\tan@ii=##2\relax
      \def\\####1####2####3{}%
     \fi}%
    \slopes@
    \slope@i=\secondy@ \advance\slope@i-\firsty@
    \ifN@\else\multiply\slope@i\m@ne\fi
    \multiply\slope@i\tan@ii \multiply\slope@i\tan@iip \multiply\slope@i\tw@
    \count@\tan@i \multiply\count@\tan@iip
    \extracount@\tan@ip \multiply\extracount@\tan@ii
    \advance\count@\extracount@
    \slope@ii=\secondx@ \advance\slope@ii-\firstx@
    \ifE@\else\multiply\slope@ii\m@ne\fi
    \multiply\slope@ii\count@
    \ifdim\slope@i<\slope@ii
     \global\tan@i=\tan@ip \global\tan@ii=\tan@iip
     \global\advance\angcount@\m@ne
    \fi
   \fi
  \fi
 \fi}%
}
\def\slope@a#1{{\def\\##1##2##3{\ifnum##3=#1\global\tan@i=##1\relax
 \global\tan@ii=##2\relax\fi}\slopes@}}
\newcount\i@
\newcount\j@
\newcount\colcount@
\newcount\Colcount@
\newcount\tcolcount@
\newdimen\rowht@
\newdimen\rowdp@
\newcount\rowcount@
\newcount\Rowcount@
\newcount\maxcolrow@
\newtoks\colwidthtoks@
\newtoks\Rowheighttoks@
\newtoks\Rowdepthtoks@
\newtoks\widthtoks@
\newtoks\Widthtoks@
\newtoks\heighttoks@
\newtoks\Heighttoks@
\newtoks\depthtoks@
\newtoks\Depthtoks@
\newif\iffirstnewCDcr@
\def\dotoks@i{%
 \global\widthtoks@=\expandafter{\the\widthtoks@\else\getdim@\z@\fi}%
 \global\heighttoks@=\expandafter{\the\heighttoks@\else\getdim@\z@\fi}%
 \global\depthtoks@=\expandafter{\the\depthtoks@\else\getdim@\z@\fi}}
\def\dotoks@ii{%
 \global\widthtoks@{\ifcase\j@}%
 \global\heighttoks@{\ifcase\j@}%
 \global\depthtoks@{\ifcase\j@}}
\def\prenewCD@#1\endnewCD{\setbox\zer@
 \vbox{%
  \def\arrow@##1##2{{}}%
  \rowcount@\m@ne \colcount@\z@ \Colcount@\z@
  \firstnewCDcr@true \toks@{}%
  \widthtoks@{\ifcase\j@}%
  \Widthtoks@{\ifcase\i@}%
  \heighttoks@{\ifcase\j@}%
  \Heighttoks@{\ifcase\i@}%
  \depthtoks@{\ifcase\j@}%
  \Depthtoks@{\ifcase\i@}%
  \Rowheighttoks@{\ifcase\i@}%
  \Rowdepthtoks@{\ifcase\i@}%
  \Let@
  \everycr{%
   \noalign{%
    \global\advance\rowcount@\@ne
    \ifnum\colcount@<\Colcount@
    \else
     \global\Colcount@=\colcount@ \global\maxcolrow@=\rowcount@
    \fi
    \global\colcount@\z@
    \iffirstnewCDcr@
     \global\firstnewCDcr@false
    \else
     \edef\next@{\the\Rowheighttoks@\noexpand\or\noexpand\getdim@\the\rowht@}%
      \global\Rowheighttoks@=\expandafter{\next@}%
     \edef\next@{\the\Rowdepthtoks@\noexpand\or\noexpand\getdim@\the\rowdp@}%
      \global\Rowdepthtoks@=\expandafter{\next@}%
     \global\rowht@\z@ \global\rowdp@\z@
     \dotoks@i
     \edef\next@{\the\Widthtoks@\noexpand\or\the\widthtoks@}%
      \global\Widthtoks@=\expandafter{\next@}%
     \edef\next@{\the\Heighttoks@\noexpand\or\the\heighttoks@}%
      \global\Heighttoks@=\expandafter{\next@}%
     \edef\next@{\the\Depthtoks@\noexpand\or\the\depthtoks@}%
      \global\Depthtoks@=\expandafter{\next@}%
     \dotoks@ii
    \fi}}%
  \tabskip\z@
  \halign{&\setbox\zer@\hbox{\vrule height10\p@ width\z@ depth\z@
   $\m@th\displaystyle{##}$}\copy\zer@
   \ifdim\ht\zer@>\rowht@ \global\rowht@\ht\zer@ \fi
   \ifdim\dp\zer@>\rowdp@ \global\rowdp@\dp\zer@ \fi
   \global\advance\colcount@\@ne
   \edef\next@{\the\widthtoks@\noexpand\or\noexpand\getdim@\the\wd\zer@}%
    \global\widthtoks@=\expandafter{\next@}%
   \edef\next@{\the\heighttoks@\noexpand\or\noexpand\getdim@\the\ht\zer@}%
    \global\heighttoks@=\expandafter{\next@}%
   \edef\next@{\the\depthtoks@\noexpand\or\noexpand\getdim@\the\dp\zer@}%
    \global\depthtoks@=\expandafter{\next@}%
   \cr#1\crcr}}%
 \Rowcount@=\rowcount@
 \global\Widthtoks@=\expandafter{\the\Widthtoks@\fi\relax}%
 \edef\Width@##1##2{\i@=##1\relax\j@=##2\relax\the\Widthtoks@}%
 \global\Heighttoks@=\expandafter{\the\Heighttoks@\fi\relax}%
 \edef\Height@##1##2{\i@=##1\relax\j@=##2\relax\the\Heighttoks@}%
 \global\Depthtoks@=\expandafter{\the\Depthtoks@\fi\relax}%
 \edef\Depth@##1##2{\i@=##1\relax\j@=##2\relax\the\Depthtoks@}%
 \edef\next@{\the\Rowheighttoks@\noexpand\fi\relax}%
 \global\Rowheighttoks@=\expandafter{\next@}%
 \edef\Rowheight@##1{\i@=##1\relax\the\Rowheighttoks@}%
 \edef\next@{\the\Rowdepthtoks@\noexpand\fi\relax}%
 \global\Rowdepthtoks@=\expandafter{\next@}%
 \edef\Rowdepth@##1{\i@=##1\relax\the\Rowdepthtoks@}%
 \colwidthtoks@{\fi}%
 \setbox\zer@\vbox{%
  \unvbox\zer@
  \count@\rowcount@
  \loop
   \unskip\unpenalty
   \setbox\zer@\lastbox
   \ifnum\count@>\maxcolrow@ \advance\count@\m@ne
   \repeat
  \hbox{%
   \unhbox\zer@
   \count@\z@
   \loop
    \unskip
    \setbox\zer@\lastbox
    \edef\next@{\noexpand\or\noexpand\getdim@\the\wd\zer@\the\colwidthtoks@}%
     \global\colwidthtoks@=\expandafter{\next@}%
    \advance\count@\@ne
    \ifnum\count@<\Colcount@
    \repeat}}%
 \edef\next@{\noexpand\ifcase\noexpand\i@\the\colwidthtoks@}%
  \global\colwidthtoks@=\expandafter{\next@}%
 \edef\Colwidth@##1{\i@=##1\relax\the\colwidthtoks@}%
 \colwidthtoks@{}\Rowheighttoks@{}\Rowdepthtoks@{}\widthtoks@{}%
 \Widthtoks@{}\heighttoks@{}\Heighttoks@{}\depthtoks@{}\Depthtoks@{}%
}
\newcount\xoff@
\newcount\yoff@
\newcount\endcount@
\newcount\rcount@
\newdimen\firstx@
\newdimen\firsty@
\newdimen\secondx@
\newdimen\secondy@
\newdimen\tocenter@
\newdimen\charht@
\newdimen\charwd@
\def\outside@{\Err@{This arrow points outside the \string\newCD}}
\newif\ifsvertex@
\newif\iftvertex@
\def\arrow@#1#2{\xoff@=#1\relax\yoff@=#2\relax
 \count@\rowcount@ \advance\count@-\yoff@
 \ifnum\count@<\@ne \outside@ \else \ifnum\count@>\Rowcount@ \outside@ \fi\fi
 \count@\colcount@ \advance\count@\xoff@
 \ifnum\count@<\@ne \outside@ \else \ifnum\count@>\Colcount@ \outside@\fi\fi
 \tcolcount@\colcount@ \advance\tcolcount@\xoff@
 \Width@\rowcount@\colcount@ \tocenter@=-\getdim@ \divide\tocenter@\tw@
 \ifdim\getdim@=\z@
  \firstx@\z@ \firsty@\mathaxis@ \svertex@true
 \else
  \svertex@false
  \ifHshort@
   \Colwidth@\colcount@
    \ifE@ \firstx@=.5\getdim@ \else \firstx@=-.5\getdim@ \fi
  \else
   \ifE@ \firstx@=\getdim@ \else \firstx@=-\getdim@ \fi
   \divide\firstx@\tw@
  \fi
  \ifE@
   \ifH@ \advance\firstx@\thr@@\p@ \else \advance\firstx@-\thr@@\p@ \fi
  \else
   \ifH@ \advance\firstx@-\thr@@\p@ \else \advance\firstx@\thr@@\p@ \fi
  \fi
  \ifN@
   \Height@\rowcount@\colcount@ \firsty@=\getdim@
   \ifV@ \advance\firsty@\thr@@\p@ \fi
  \else
   \ifV@
    \Depth@\rowcount@\colcount@ \firsty@=-\getdim@
    \advance\firsty@-\thr@@\p@
   \else
    \firsty@\z@
   \fi
  \fi
 \fi
 \ifV@
 \else
  \Colwidth@\colcount@
  \ifE@ \secondx@=\getdim@ \else \secondx@=-\getdim@ \fi
  \divide\secondx@\tw@
  \ifE@ \else \getcgap@\colcount@ \advance\secondx@-\getdim@ \fi
  \endcount@=\colcount@ \advance\endcount@\xoff@
  \count@=\colcount@
  \ifE@
   \advance\count@\@ne
   \loop
    \ifnum\count@<\endcount@
    \Colwidth@\count@ \advance\secondx@\getdim@
    \getcgap@\count@ \advance\secondx@\getdim@
    \advance\count@\@ne
    \repeat
  \else
   \advance\count@\m@ne
   \loop
    \ifnum\count@>\endcount@
    \Colwidth@\count@ \advance\secondx@-\getdim@
    \getcgap@\count@ \advance\secondx@-\getdim@
    \advance\count@\m@ne
    \repeat
  \fi
  \Colwidth@\count@ \divide\getdim@\tw@
  \ifHshort@
  \else
   \ifE@ \advance\secondx@\getdim@ \else \advance\secondx@-\getdim@ \fi
  \fi
  \ifE@ \getcgap@\count@ \advance\secondx@\getdim@ \fi
  \rcount@\rowcount@ \advance\rcount@-\yoff@
  \Width@\rcount@\count@ \divide\getdim@\tw@
  \tvertex@false
  \ifH@\ifdim\getdim@=\z@\tvertex@true\Hshort@false\fi\fi
  \ifHshort@
  \else
   \ifE@ \advance\secondx@-\getdim@ \else \advance\secondx@\getdim@ \fi
  \fi
  \iftvertex@
   \advance\secondx@.4\p@
  \else
   \ifE@ \advance\secondx@-\thr@@\p@ \else \advance\secondx@\thr@@\p@ \fi
  \fi
 \fi
 \ifH@
 \else
  \ifN@
   \Rowheight@\rowcount@ \secondy@\getdim@
  \else
   \Rowdepth@\rowcount@ \secondy@-\getdim@
   \getrgap@\rowcount@ \advance\secondy@-\getdim@
  \fi
  \endcount@=\rowcount@ \advance\endcount@-\yoff@
  \count@=\rowcount@
  \ifN@
   \advance\count@\m@ne
   \loop
    \ifnum\count@>\endcount@
    \Rowheight@\count@ \advance\secondy@\getdim@
    \Rowdepth@\count@ \advance\secondy@\getdim@
    \getrgap@\count@ \advance\secondy@\getdim@
    \advance\count@\m@ne
    \repeat
  \else
   \advance\count@\@ne
   \loop
    \ifnum\count@<\endcount@
    \Rowheight@\count@ \advance\secondy@-\getdim@
    \Rowdepth@\count@ \advance\secondy@-\getdim@
    \getrgap@\count@ \advance\secondy@-\getdim@
    \advance\count@\@ne
    \repeat
  \fi
  \tvertex@false
  \ifV@\Width@\count@\colcount@\ifdim\getdim@=\z@\tvertex@true\fi\fi
  \ifN@
   \getrgap@\count@ \advance\secondy@\getdim@
   \Rowdepth@\count@ \advance\secondy@\getdim@
   \iftvertex@
    \advance\secondy@\mathaxis@
   \else
    \Depth@\count@\tcolcount@ \advance\secondy@-\getdim@
    \advance\secondy@-\thr@@\p@
   \fi
  \else
   \Rowheight@\count@ \advance\secondy@-\getdim@
   \iftvertex@
    \advance\secondy@\mathaxis@
   \else
    \Height@\count@\tcolcount@ \advance\secondy@\getdim@
    \advance\secondy@\thr@@\p@
   \fi
  \fi
 \fi
 \ifV@\else\advance\firstx@\sxdimen@\fi
 \ifH@\else\advance\firsty@\sydimen@\fi
 \iftX@
  \advance\secondy@\tXdimen@ii
  \advance\secondx@\tXdimen@i
  \slope@
 \else
  \iftY@
   \advance\secondy@\tYdimen@ii
   \advance\secondx@\tYdimen@i
   \slope@
   \secondy@=\secondx@ \advance\secondy@-\firstx@
   \ifNESW@ \else \multiply\secondy@\m@ne \fi
   \multiply\secondy@\tan@i \divide\secondy@\tan@ii \advance\secondy@\firsty@
  \else
   \ifa@
    \slope@
    \ifNESW@ \global\advance\angcount@\exacount@ \else
      \global\advance\angcount@-\exacount@ \fi
    \ifnum\angcount@>23 \angcount@23 \fi
    \ifnum\angcount@<\@ne \angcount@\@ne \fi
    \slope@a\angcount@
    \ifY@
     \advance\secondy@\Ydimen@
    \else
     \ifX@
      \advance\secondx@\Xdimen@
      \dimen@\secondx@ \advance\dimen@-\firstx@
      \ifNESW@\else\multiply\dimen@\m@ne\fi
      \multiply\dimen@\tan@i \divide\dimen@\tan@ii
      \advance\dimen@\firsty@ \secondy@=\dimen@
     \fi
    \fi
   \else
    \ifH@\else\ifV@\else\slope@\fi\fi
   \fi
  \fi
 \fi
 \ifH@\else\ifV@\else\ifsvertex@\else
  \dimen@=6\p@ \multiply\dimen@\tan@ii
  \count@=\tan@i \advance\count@\tan@ii \divide\dimen@\count@
  \ifE@ \advance\firstx@\dimen@ \else \advance\firstx@-\dimen@ \fi
  \multiply\dimen@\tan@i \divide\dimen@\tan@ii
  \ifN@ \advance\firsty@\dimen@ \else \advance\firsty@-\dimen@ \fi
 \fi\fi\fi
 \ifp@
  \ifH@\else\ifV@\else
   \getcos@\pdimen@ \advance\firsty@\dimen@ \advance\secondy@\dimen@
   \ifNESW@ \advance\firstx@-\dimen@ii \else \advance\firstx@\dimen@ii \fi
  \fi\fi
 \fi
 \ifH@\else\ifV@\else
  \ifnum\tan@i>\tan@ii
   \charht@=10\p@ \charwd@=10\p@
   \multiply\charwd@\tan@ii \divide\charwd@\tan@i
  \else
   \charwd@=10\p@ \charht@=10\p@
   \divide\charht@\tan@ii \multiply\charht@\tan@i
  \fi
  \ifnum\tcount@=\thr@@
   \ifN@ \advance\secondy@-.3\charht@ \else\advance\secondy@.3\charht@ \fi
  \fi
  \ifnum\scount@=\tw@
   \ifE@ \advance\firstx@.3\charht@ \else \advance\firstx@-.3\charht@ \fi
  \fi
  \ifnum\tcount@=12
   \ifN@ \advance\secondy@-\charht@ \else \advance\secondy@\charht@ \fi
  \fi
  \iftY@
  \else
   \ifa@
    \ifX@
    \else
     \secondx@\secondy@ \advance\secondx@-\firsty@
     \ifNESW@\else\multiply\secondx@\m@ne\fi
     \multiply\secondx@\tan@ii \divide\secondx@\tan@i
     \advance\secondx@\firstx@
    \fi
   \fi
  \fi
 \fi\fi
 \ifH@\harrow@\else\ifV@\varrow@\else\arrow@@\fi\fi}
\newdimen\mathaxis@
\mathaxis@90\p@ \divide\mathaxis@36
\def\harrow@b{\ifE@\hskip\tocenter@\hskip\firstx@\fi}
\def\harrow@bb{\ifE@\hskip\xdimen@\else\hskip\Xdimen@\fi}
\def\harrow@e{\ifE@\else\hskip-\firstx@\hskip-\tocenter@\fi}
\def\harrow@ee{\ifE@\hskip-\Xdimen@\else\hskip-\xdimen@\fi}
\def\harrow@{\dimen@\secondx@\advance\dimen@-\firstx@
 \ifE@ \let\next@\rlap \else  \multiply\dimen@\m@ne \let\next@\llap \fi
 \next@{%
  \harrow@b
  \smash{\raise\pdimen@\hbox to\dimen@
   {\harrow@bb\arrow@ii
    \ifnum\arrcount@=\m@ne \else \ifnum\arrcount@=\thr@@ \else
     \ifE@
      \ifnum\scount@=\m@ne
      \else
       \ifcase\scount@\or\or\char118 \or\char117 \or\or\or\char119 \or
       \char120 \or\char121 \or\char122 \or\or\or\arrow@i\char125 \or
       \char117 \hskip\thr@@\p@\char117 \hskip-\thr@@\p@\fi
      \fi
     \else
      \ifnum\tcount@=\m@ne
      \else
       \ifcase\tcount@\char117 \or\or\char117 \or\char118 \or\char119 \or
       \char120\or\or\or\or\or\char121 \or\char122 \or\arrow@i\char125
       \or\char117 \hskip\thr@@\p@\char117 \hskip-\thr@@\p@\fi
      \fi
     \fi
    \fi\fi
    \dimen@\mathaxis@ \advance\dimen@.2\p@
    \dimen@ii\mathaxis@ \advance\dimen@ii-.2\p@
    \ifnum\arrcount@=\m@ne
     \let\leads@\null
    \else
     \ifcase\arrcount@
      \def\leads@{\hrule height\dimen@ depth-\dimen@ii}\or
      \def\leads@{\hrule height\dimen@ depth-\dimen@ii}\or
      \def\leads@{\hbox to10\p@{%
       \leaders\hrule height\dimen@ depth-\dimen@ii\hfil
       \hfil
      \leaders\hrule height\dimen@ depth-\dimen@ii\hskip\z@ plus2fil\relax
       \hfil
       \leaders\hrule height\dimen@ depth-\dimen@ii\hfil}}\or
     \def\leads@{\hbox{\hbox to10\p@{\dimen@\mathaxis@ \advance\dimen@1.2\p@
       \dimen@ii\dimen@ \advance\dimen@ii-.4\p@
       \leaders\hrule height\dimen@ depth-\dimen@ii\hfil}%
       \kern-10\p@
       \hbox to10\p@{\dimen@\mathaxis@ \advance\dimen@-1.2\p@
       \dimen@ii\dimen@ \advance\dimen@ii-.4\p@
       \leaders\hrule height\dimen@ depth-\dimen@ii\hfil}}}\fi
    \fi
    \cleaders\leads@\hfil
    \ifnum\arrcount@=\m@ne\else\ifnum\arrcount@=\thr@@\else
     \arrow@i
     \ifE@
      \ifnum\tcount@=\m@ne
      \else
       \ifcase\tcount@\char119 \or\or\char119 \or\char120 \or\char121 \or
       \char122 \or \or\or\or\or\char123\or\char124 \or
       \char125 \or\char119 \hskip-\thr@@\p@\char119 \hskip\thr@@\p@\fi
      \fi
     \else
      \ifcase\scount@\or\or\char120 \or\char119 \or\or\or\char121 \or\char122
      \or\char123 \or\char124 \or\or\or\char125 \or
      \char119 \hskip-\thr@@\p@\char119 \hskip\thr@@\p@\fi
     \fi
    \fi\fi
    \harrow@ee}}%
  \harrow@e}%
 \iflabel@i
  \dimen@ii\z@ \setbox\zer@\hbox{$\m@th\tsize@@\label@i$}%
  \ifnum\arrcount@=\m@ne
  \else
   \advance\dimen@ii\mathaxis@
   \advance\dimen@ii\dp\zer@ \advance\dimen@ii\tw@\p@
   \ifnum\arrcount@=\thr@@ \advance\dimen@ii\tw@\p@ \fi
  \fi
  \advance\dimen@ii\pdimen@
  \next@{\harrow@b\smash{\raise\dimen@ii\hbox to\dimen@
   {\harrow@bb\hskip\tw@\ldimen@i\hfil\box\zer@\hfil\harrow@ee}}\harrow@e}%
 \fi
 \iflabel@ii
  \ifnum\arrcount@=\m@ne
  \else
   \setbox\zer@\hbox{$\m@th\tsize@\label@ii$}%
   \dimen@ii-\ht\zer@ \advance\dimen@ii-\tw@\p@
   \ifnum\arrcount@=\thr@@ \advance\dimen@ii-\tw@\p@ \fi
   \advance\dimen@ii\mathaxis@ \advance\dimen@ii\pdimen@
   \next@{\harrow@b\smash{\raise\dimen@ii\hbox to\dimen@
    {\harrow@bb\hskip\tw@\ldimen@ii\hfil\box\zer@\hfil\harrow@ee}}\harrow@e}%
  \fi
 \fi}
\let\tsize@\tsize
\def\tsizenewCDlabels{\let\tsize@\tsize}
\def\ssizenewCDlabels{\let\tsize@\ssize}
\def\tsize@@{\ifnum\arrcount@=\m@ne\else\tsize@\fi}
\def\varrow@{\dimen@\secondy@ \advance\dimen@-\firsty@
 \ifN@ \else \multiply\dimen@\m@ne \fi
 \setbox\zer@\vbox to\dimen@
  {\ifN@ \vskip-\Ydimen@ \else \vskip\ydimen@ \fi
   \ifnum\arrcount@=\m@ne\else\ifnum\arrcount@=\thr@@\else
    \hbox{\arrow@iii
     \ifN@
      \ifnum\tcount@=\m@ne
      \else
       \ifcase\tcount@\char117 \or\or\char117 \or\char118 \or\char119 \or
       \char120 \or\or\or\or\or\char121 \or\char122 \or\char123 \or
       \vbox{\hbox{\char117 }\nointerlineskip\vskip\thr@@\p@
       \hbox{\char117 }\vskip-\thr@@\p@}\fi
      \fi
     \else
      \ifcase\scount@\or\or\char118 \or\char117 \or\or\or\char119 \or
      \char120 \or\char121 \or\char122 \or\or\or\char123 \or
      \vbox{\hbox{\char117 }\nointerlineskip\vskip\thr@@\p@
      \hbox{\char117 }\vskip-\thr@@\p@}\fi
     \fi}%
    \nointerlineskip
   \fi\fi
   \ifnum\arrcount@=\m@ne
    \let\leads@\null
   \else
    \ifcase\arrcount@\let\leads@\vrule\or\let\leads@\vrule\or
    \def\leads@{\vbox to10\p@{%
     \hrule height 1.67\p@ depth\z@ width.4\p@
     \vfil
     \hrule height 3.33\p@ depth\z@ width.4\p@
     \vfil
     \hrule height 1.67\p@ depth\z@ width.4\p@}}\or
    \def\leads@{\hbox{\vrule height\p@\hskip\tw@\p@\vrule}}\fi
   \fi
  \cleaders\leads@\vfill\nointerlineskip
   \ifnum\arrcount@=\m@ne\else\ifnum\arrcount@=\thr@@\else
    \hbox{\arrow@iv
     \ifN@
      \ifcase\scount@\or\or\char118 \or\char117 \or\or\or\char119 \or
      \char120 \or\char121 \or\char122 \or\or\or\arrow@iii\char123 \or
      \vbox{\hbox{\char117 }\nointerlineskip\vskip-\thr@@\p@
      \hbox{\char117 }\vskip\thr@@\p@}\fi
     \else
      \ifnum\tcount@=\m@ne
      \else
       \ifcase\tcount@\char117 \or\or\char117 \or\char118 \or\char119 \or
       \char120 \or\or\or\or\or\char121 \or\char122 \or\arrow@iii\char123 \or
       \vbox{\hbox{\char117 }\nointerlineskip\vskip-\thr@@\p@
       \hbox{\char117 }\vskip\thr@@\p@}\fi
      \fi
     \fi}%
   \fi\fi
   \ifN@\vskip\ydimen@\else\vskip-\Ydimen@\fi}%
 \ifN@
  \dimen@ii\firsty@
 \else
  \dimen@ii-\firsty@ \advance\dimen@ii\ht\zer@ \multiply\dimen@ii\m@ne
 \fi
 \rlap{\smash{\hskip\tocenter@ \hskip\pdimen@ \raise\dimen@ii \box\zer@}}%
 \iflabel@i
  \setbox\zer@\vbox to\dimen@{\vfil
   \hbox{$\m@th\tsize@@\label@i$}\vskip\tw@\ldimen@i\vfil}%
  \rlap{\smash{\hskip\tocenter@ \hskip\pdimen@
  \ifnum\arrcount@=\m@ne \let\next@\relax \else \let\next@\llap \fi
  \next@{\raise\dimen@ii\hbox{\ifnum\arrcount@=\m@ne \hskip-.5\wd\zer@ \fi
   \box\zer@ \ifnum\arrcount@=\m@ne \else \hskip\tw@\p@ \fi}}}}%
 \fi
 \iflabel@ii
  \ifnum\arrcount@=\m@ne
  \else
   \setbox\zer@\vbox to\dimen@{\vfil
    \hbox{$\m@th\tsize@\label@ii$}\vskip\tw@\ldimen@ii\vfil}%
   \rlap{\smash{\hskip\tocenter@ \hskip\pdimen@
   \rlap{\raise\dimen@ii\hbox{\ifnum\arrcount@=\thr@@ \hskip4.5\p@ \else
    \hskip2.5\p@ \fi\box\zer@}}}}%
  \fi
 \fi
}
\newdimen\goal@
\newdimen\shifted@
\newcount\Tcount@
\newcount\Scount@
\newbox\shaft@
\newcount\slcount@
\def\getcos@#1{%
 \ifnum\tan@i<\tan@ii
  \dimen@#1%
  \ifnum\slcount@<8 \count@9 \else \ifnum\slcount@<12 \count@8 \else
   \count@7 \fi\fi
  \multiply\dimen@\count@ \divide\dimen@10
  \dimen@ii\dimen@ \multiply\dimen@ii\tan@i \divide\dimen@ii\tan@ii
 \else
  \dimen@ii#1%
  \count@-\slcount@ \advance\count@24
  \ifnum\count@<8 \count@9 \else \ifnum\count@<12 \count@8
   \else\count@7 \fi\fi
  \multiply\dimen@ii\count@ \divide\dimen@ii10
  \dimen@\dimen@ii \multiply\dimen@\tan@ii \divide\dimen@\tan@i
 \fi}
\newdimen\adjust@
\def\Nnext@{\ifN@\let\next@\raise\else\let\next@\lower\fi}
\def\arrow@@{\slcount@\angcount@
 \ifNESW@
  \ifnum\angcount@<10
   \let\arrowfont@=\arrow@i \advance\angcount@\m@ne \multiply\angcount@13
  \else
   \ifnum\angcount@<19
    \let\arrowfont@=\arrow@ii \advance\angcount@-10 \multiply\angcount@13
   \else
    \let\arrowfont@=\arrow@iii \advance\angcount@-19 \multiply\angcount@13
  \fi\fi
  \Tcount@\angcount@
 \else
  \ifnum\angcount@<5
   \let\arrowfont@=\arrow@iii \advance\angcount@\m@ne \multiply\angcount@13
   \advance\angcount@65
  \else
   \ifnum\angcount@<14
    \let\arrowfont@=\arrow@iv \advance\angcount@-5 \multiply\angcount@13
   \else
    \ifnum\angcount@<23
     \let\arrowfont@=\arrow@v \advance\angcount@-14 \multiply\angcount@13
    \else
     \let\arrowfont@=\arrow@i \angcount@=117
  \fi\fi\fi
  \ifnum\angcount@=117 \Tcount@=115 \else\Tcount@\angcount@ \fi
 \fi
 \Scount@\Tcount@
 \ifE@
  \ifnum\tcount@=\z@ \advance\Tcount@\tw@ \else\ifnum\tcount@=13
   \advance\Tcount@\tw@ \else \advance\Tcount@\tcount@ \fi\fi
  \ifnum\scount@=\z@ \else \ifnum\scount@=13 \advance\Scount@\thr@@ \else
   \advance\Scount@\scount@ \fi\fi
 \else
  \ifcase\tcount@\advance\Tcount@\thr@@\or\or\advance\Tcount@\thr@@\or
  \advance\Tcount@\tw@\or\advance\Tcount@6 \or\advance\Tcount@7
  \or\or\or\or\or \advance\Tcount@8 \or\advance\Tcount@9 \or
  \advance\Tcount@12 \or\advance\Tcount@\thr@@\fi
  \ifcase\scount@\or\or\advance\Scount@\thr@@\or\advance\Scount@\tw@\or
  \or\or\advance\Scount@4 \or\advance\Scount@5 \or\advance\Scount@10
  \or\advance\Scount@11 \or\or\or\advance\Scount@12 \or\advance
  \Scount@\tw@\fi
 \fi
 \ifcase\arrcount@\or\or\advance\angcount@\@ne\else\fi
 \ifN@ \shifted@=\firsty@ \else\shifted@=-\firsty@ \fi
 \ifE@ \else\advance\shifted@\charht@ \fi
 \goal@=\secondy@ \advance\goal@-\firsty@
 \ifN@\else\multiply\goal@\m@ne\fi
 \setbox\shaft@\hbox{\arrowfont@\char\angcount@}%
 \ifnum\arrcount@=\thr@@
  \getcos@{1.5\p@}%
  \setbox\shaft@\hbox to\wd\shaft@{\arrowfont@
   \rlap{\hskip\dimen@ii
    \smash{\ifNESW@\let\next@\lower\else\let\next@\raise\fi
     \next@\dimen@\hbox{\arrowfont@\char\angcount@}}}%
   \rlap{\hskip-\dimen@ii
    \smash{\ifNESW@\let\next@\raise\else\let\next@\lower\fi
      \next@\dimen@\hbox{\arrowfont@\char\angcount@}}}\hfil}%
 \fi
 \rlap{\smash{\hskip\tocenter@\hskip\firstx@
  \ifnum\arrcount@=\m@ne
  \else
   \ifnum\arrcount@=\thr@@
   \else
    \ifnum\scount@=\m@ne
    \else
     \ifnum\scount@=\z@
     \else
      \setbox\zer@\hbox{\ifnum\angcount@=117 \arrow@v\else\arrowfont@\fi
       \char\Scount@}%
      \ifNESW@
       \ifnum\scount@=\tw@
        \dimen@=\shifted@ \advance\dimen@-\charht@
        \ifN@\hskip-\wd\zer@\fi
        \Nnext@
        \next@\dimen@\copy\zer@
        \ifN@\else\hskip-\wd\zer@\fi
       \else
        \Nnext@
        \ifN@\else\hskip-\wd\zer@\fi
        \next@\shifted@\copy\zer@
        \ifN@\hskip-\wd\zer@\fi
       \fi
       \ifnum\scount@=12
        \advance\shifted@\charht@ \advance\goal@-\charht@
        \ifN@ \hskip\wd\zer@ \else \hskip-\wd\zer@ \fi
       \fi
       \ifnum\scount@=13
        \getcos@{\thr@@\p@}%
        \ifN@ \hskip\dimen@ \else \hskip-\wd\zer@ \hskip-\dimen@ \fi
        \adjust@\shifted@ \advance\adjust@\dimen@ii
        \Nnext@
        \next@\adjust@\copy\zer@
        \ifN@ \hskip-\dimen@ \hskip-\wd\zer@ \else \hskip\dimen@ \fi
       \fi
      \else
       \ifN@\hskip-\wd\zer@\fi
       \ifnum\scount@=\tw@
        \ifN@ \hskip\wd\zer@ \else \hskip-\wd\zer@ \fi
        \dimen@=\shifted@ \advance\dimen@-\charht@
        \Nnext@
        \next@\dimen@\copy\zer@
        \ifN@\hskip-\wd\zer@\fi
       \else
        \Nnext@
        \next@\shifted@\copy\zer@
        \ifN@\else\hskip-\wd\zer@\fi
       \fi
       \ifnum\scount@=12
        \advance\shifted@\charht@ \advance\goal@-\charht@
        \ifN@ \hskip-\wd\zer@ \else \hskip\wd\zer@ \fi
       \fi
       \ifnum\scount@=13
        \getcos@{\thr@@\p@}%
        \ifN@ \hskip-\wd\zer@ \hskip-\dimen@ \else \hskip\dimen@ \fi
        \adjust@\shifted@ \advance\adjust@\dimen@ii
        \Nnext@
        \next@\adjust@\copy\zer@
        \ifN@ \hskip\dimen@ \else \hskip-\dimen@ \hskip-\wd\zer@ \fi
       \fi      
      \fi
  \fi\fi\fi\fi
  \ifnum\arrcount@=\m@ne
  \else
   \loop
    \ifdim\goal@>\charht@
    \ifE@\else\hskip-\charwd@\fi
    \Nnext@
    \next@\shifted@\copy\shaft@
    \ifE@\else\hskip-\charwd@\fi
    \advance\shifted@\charht@ \advance\goal@ -\charht@
    \repeat
   \ifdim\goal@>\z@
    \dimen@=\charht@ \advance\dimen@-\goal@
    \divide\dimen@\tan@i \multiply\dimen@\tan@ii
    \ifE@ \hskip-\dimen@ \else \hskip-\charwd@ \hskip\dimen@ \fi
    \adjust@=\shifted@ \advance\adjust@-\charht@ \advance\adjust@\goal@
    \Nnext@
    \next@\adjust@\copy\shaft@
    \ifE@ \else \hskip-\charwd@ \fi
   \else
    \adjust@=\shifted@ \advance\adjust@-\charht@
   \fi
  \fi
  \ifnum\arrcount@=\m@ne
  \else
   \ifnum\arrcount@=\thr@@
   \else
    \ifnum\tcount@=\m@ne
    \else
     \setbox\zer@
      \hbox{\ifnum\angcount@=117 \arrow@v\else\arrowfont@\fi\char\Tcount@}%
     \ifnum\tcount@=\thr@@
      \advance\adjust@\charht@
      \ifE@\else\ifN@\hskip-\charwd@\else\hskip-\wd\zer@\fi\fi
     \else
      \ifnum\tcount@=12
       \advance\adjust@\charht@
       \ifE@\else\ifN@\hskip-\charwd@\else\hskip-\wd\zer@\fi\fi
      \else
       \ifE@\hskip-\wd\zer@\fi
     \fi\fi
     \Nnext@
     \next@\adjust@\copy\zer@
     \ifnum\tcount@=13
      \hskip-\wd\zer@
      \getcos@{\thr@@\p@}%
      \ifE@\hskip-\dimen@ \else\hskip\dimen@ \fi
      \advance\adjust@-\dimen@ii
      \Nnext@
      \next@\adjust@\box\zer@
     \fi
  \fi\fi\fi}}%
 \iflabel@i
  \rlap{\hskip\tocenter@
  \dimen@\firstx@ \advance\dimen@\secondx@ \divide\dimen@\tw@
  \advance\dimen@\ldimen@i
  \dimen@ii\firsty@ \advance\dimen@ii\secondy@ \divide\dimen@ii\tw@
  \multiply\ldimen@i\tan@i \divide\ldimen@i\tan@ii
  \ifNESW@ \advance\dimen@ii\ldimen@i \else \advance\dimen@ii-\ldimen@i \fi
  \setbox\zer@\hbox{\ifNESW@\else\ifnum\arrcount@=\thr@@\hskip4\p@\else
   \hskip\tw@\p@\fi\fi
   $\m@th\tsize@@\label@i$\ifNESW@\ifnum\arrcount@=\thr@@\hskip4\p@\else
   \hskip\tw@\p@\fi\fi}%
  \ifnum\arrcount@=\m@ne
   \ifNESW@ \advance\dimen@.5\wd\zer@ \advance\dimen@\p@ \else
    \advance\dimen@-.5\wd\zer@ \advance\dimen@-\p@ \fi
   \advance\dimen@ii-.5\ht\zer@
  \else
   \advance\dimen@ii\dp\zer@
   \ifnum\slcount@<6 \advance\dimen@ii\tw@\p@ \fi
  \fi
  \hskip\dimen@
  \ifNESW@ \let\next@\llap \else\let\next@\rlap \fi
  \next@{\smash{\raise\dimen@ii\box\zer@}}}%
 \fi
 \iflabel@ii
  \ifnum\arrcount@=\m@ne
  \else
   \rlap{\hskip\tocenter@
   \dimen@\firstx@ \advance\dimen@\secondx@ \divide\dimen@\tw@
   \ifNESW@ \advance\dimen@\ldimen@ii \else \advance\dimen@-\ldimen@ii \fi
   \dimen@ii\firsty@ \advance\dimen@ii\secondy@ \divide\dimen@ii\tw@
   \multiply\ldimen@ii\tan@i \divide\ldimen@ii\tan@ii
   \advance\dimen@ii\ldimen@ii
   \setbox\zer@\hbox{\ifNESW@\ifnum\arrcount@=\thr@@\hskip4\p@\else
    \hskip\tw@\p@\fi\fi
    $\m@th\tsize@\label@ii$\ifNESW@\else\ifnum\arrcount@=\thr@@\hskip4\p@
    \else\hskip\tw@\p@\fi\fi}%
   \advance\dimen@ii-\ht\zer@
   \ifnum\slcount@<9 \advance\dimen@ii-\thr@@\p@ \fi
   \ifNESW@ \let\next@\rlap \else \let\next@\llap \fi
   \hskip\dimen@\next@{\smash{\raise\dimen@ii\box\zer@}}}%
  \fi
 \fi
}
\def\outnewCD@#1{\def#1{\Err@{\string#1 must not be used within \string\newCD}}}
\newskip\prenewCDskip@
\newskip\postnewCDskip@
\prenewCDskip@\z@
\postnewCDskip@\z@
\def\prenewCDspace#1{\RIfMIfI@
 \onlydmatherr@\prenewCDspace\else\advance\prenewCDskip@#1\relax\fi\else
 \onlydmatherr@\prenewCDspace\fi}
\def\postnewCDspace#1{\RIfMIfI@
 \onlydmatherr@\postnewCDspace\else\advance\postnewCDskip@#1\relax\fi\else
 \onlydmatherr@\postnewCDspace\fi}
\def\predisplayspace#1{\RIfMIfI@
 \onlydmatherr@\predisplayspace\else
 \advance\abovedisplayskip#1\relax
 \advance\abovedisplayshortskip#1\relax\fi
 \else\onlydmatherr@\prenewCDspace\fi}
\def\postdisplayspace#1{\RIfMIfI@
 \onlydmatherr@\postdisplayspace\else
 \advance\belowdisplayskip#1\relax
 \advance\belowdisplayshortskip#1\relax\fi
 \else\onlydmatherr@\postdisplayspace\fi}
\def\PrenewCDSpace#1{\global\prenewCDskip@#1\relax}
\def\PostnewCDSpace#1{\global\postnewCDskip@#1\relax}
\def\newCD#1\endnewCD{%
 \outnewCD@\cgaps\outnewCD@\rgaps\outnewCD@\Cgaps\outnewCD@\Rgaps
 \prenewCD@#1\endnewCD
 \advance\abovedisplayskip\prenewCDskip@
 \advance\abovedisplayshortskip\prenewCDskip@
 \advance\belowdisplayskip\postnewCDskip@
 \advance\belowdisplayshortskip\postnewCDskip@
 \vcenter{\vskip\prenewCDskip@ \Let@ \colcount@\@ne \rowcount@\z@
  \everycr{%
   \noalign{%
    \ifnum\rowcount@=\Rowcount@
    \else
     \global\nointerlineskip
     \getrgap@\rowcount@ \vskip\getdim@
     \global\advance\rowcount@\@ne \global\colcount@\@ne
    \fi}}%
  \tabskip\z@
  \halign{&\global\xoff@\z@ \global\yoff@\z@
   \getcgap@\colcount@ \hskip\getdim@
   \hfil\vrule height10\p@ width\z@ depth\z@
   $\m@th\displaystyle{##}$\hfil
   \global\advance\colcount@\@ne\cr
   #1\crcr}\vskip\postnewCDskip@}%
 \prenewCDskip@\z@\postnewCDskip@\z@
 \def\getcgap@##1{\ifcase##1\or\getdim@\z@\else\getdim@\standardcgap\fi}%
 \def\getrgap@##1{\ifcase##1\getdim@\z@\else\getdim@\standardrgap\fi}%
 \let\Width@\relax\let\Height@\relax\let\Depth@\relax\let\Rowheight@\relax
 \let\Rowdepth@\relax\let\Colwdith@\relax
}
\catcode`\@=\active
\hsize 30pc
\vsize 47pc
\magnification=\magstep1
\let\[\lceil
\let\]\rfloor
\def\nmb#1#2{#2}         
\def\cit#1#2{\ifx#1!\cite{#2}\else#2\fi} 
\def\idx{}               
\def\ign#1{}             
\redefine\o{\circ}
\let\ceylon\colon
\redefine\colon{\hskip.05em\ceylon}
\define\({\big(}
\define\){\big)}

\define\al{\alpha}
\define\be{\beta}

\define\de{\delta}

\define\la{\lambda}
\define\rh{\rho}
\define\si{\sigma}

\define\ph{\varphi}

\define\ps{\psi}

\define\x{\times}
\define\Id{\operatorname{Id}}
\define\g{\frak g}
\define\h{\frak h}
\define\e{\frak e}

\redefine\l{\frak l}
\define\ad{\operatorname{ad}}
\define\pr{\operatorname{pr}}

\define\der{\operatorname{der}}

\redefine\L{\operatorname{\Cal L}}
\topmatter
\title  Extensions of Lie algebras 
\endtitle
\author Dmitri Alekseevsky, Peter W. Michor, Wolfgang Ruppert  
\endauthor
\affil
Erwin Schr\"odinger Institut f\"ur Mathematische Physik,\endgraf
Boltzmanngasse 9, A-1090 Wien, Austria
\endaffil
\address 
D.V. Alekseevsky: 
Department of Mathematics,
University of Hull,
Cottingham Road,
Hull, HU6 7RX,
England
\endaddress
\email
d.v.alekseevsky\@maths.hull.ac.uk 
\endemail
\address
P\. Michor: Institut f\"ur Mathematik, Universit\"at Wien,
Strudlhofgasse 4, A-1090 Wien, Austria
\endaddress
\email Peter.Michor\@esi.ac.at \endemail
\address
W\. Ruppert: 
Institut f\"ur Mathematik und angewandte Statistik, Universit\"at 
f\"ur Bodenkultur,
Gregor Mendelstrasse 33, A-1180 Wien, Austria
\endaddress
\email Ruppert\@edv1.boku.ac.at \endemail

\dedicatory \enddedicatory
\date {October 25, 2000} \enddate
\thanks P.W.M. was supported  
by `Fonds zur F\"orderung der wissenschaftlichen  
Forschung, Projekt P~14195~MAT'. 
\endthanks
\abstract We review (non-abelian) extensions of a given Lie algebra, 
identify a 3-dimensional cohomological obstruction to the existence of
extensions.
A striking analogy to the setting of covariant exterior derivatives, curvature, 
and the Bianchi identity in differential geometry is spelled out. 
\endabstract
\endtopmatter

\long\def\alert#1{\relax}
\document


\subhead\nmb.{1}. Introduction \endsubhead
The theory of group extensions and their interpretation in terms of 
cohomology is  well known, see 
\cit!{2}, \cit!{5}, \cit!{3}, \cit!{1}, e.g.
When we wrote this paper we thought that the counterpart for Lie 
algebras (for non-abelian extensions) was not spelled out in  
detail in the literature. We presented it here in this short note,
with special emphasis to connections with 
the (algebraic) theory of covariant exterior 
derivatives, curvature and the Bianchi identity in differential 
geometry (see section \nmb!{3}). 
Kirill Mackenzie pointed out to us, that most of our results are 
available in \cit!{4}, \cit!{12}, \cit!{15}, and generalizations 
for Lie algebroids are in \cit!{11}. 
So this paper will not appear in print. 

\subhead\nmb.{2}. Describing extensions \endsubhead
Consider any exact sequence of homomorphisms of Lie algebras: 
$$
0\to \h @>i>> \e @>p>> \g \to 0.
$$
Consider a linear mapping $s:\g \to \e$ with $p\o s=\Id_\g$. 
Then $s$ induces mappings 
$$\align 
\al&:\g\to \der(\h),\qquad
     \al_X(H)= [s(X),H], \tag{\nmb!{2}.1}\\
\rh&:\bigwedge^2\g \to \h,\qquad
     \rh(X,Y)= [s(X),s(Y)]-s([X,Y]), \tag{\nmb!{2}.2}\\
\endalign$$
which are easily seen to satisfy 
$$\align
&[\al_X,\al_Y]-\al_{[X,Y]} = \ad_{\rh(X,Y)} \tag{\nmb!{2}.3}\\
&\sum_{\text{cyclic}\{X,Y,Z\}} 
     \Bigl(\al_X \rh(Y,Z) - \rh([X,Y],Z) \Bigr). =0\tag{\nmb!{2}.4}
\endalign$$

We can completely describe the Lie algebra structure on 
$\e=\h\oplus s(\g)$ in terms of $\al$ and $\rh$ :
$$\multline
[H_1+s(X_1),H_2+s(X_2)] =\\
= ([H_1,H_2]+\al_{X_1}H_2 -\al_{X_2}H_1 
+\rh(X_1,X_2)) + s[X_1,X_2] 
\endmultline\tag{\nmb!{2}.5}$$
and one can check that formula \thetag{\nmb!{2}.5} 
gives a Lie algebra structure 
on $\h\oplus s(\g)$, if $\al:\g\to\der(\h)$ and 
$\rh:\bigwedge^2\g\to \h$ satisfy \thetag{\nmb!{2}.3} and 
\thetag{\nmb!{2}.4}. 

\subhead\nmb.{3}. Motivation: Lie algebra extensions 
associated to a principal bundle
\endsubhead

Let $\pi : P \to M = P/K $ be a principal bundle with structure 
group $K$; i.e\. $P$ is a manifold with a free right action of a Lie 
group $K$ and $\pi$ is the projection on the orbit space $M=P/K$. 
Denote by $\frak g = \frak X(M)$ the Lie algebra of the vector fields 
on $M$, by $\frak e = \frak X(P)^K$  the Lie algebra of $K$-invariant 
vector fields on $P$ and by $\frak h = \frak X_{\text{vert}}(P)^K$ the 
ideal of the $K$-invariant vertical vector fields of $\frak e$. 
Geometrically,  
$\frak e$ is the Lie algebra of infinitesimal automorphisms of the 
principal bundle $P$ 
and $\frak h$ is the ideal of infinitesimal automorphisms acting 
trivially on $M$, i.e\. the Lie algebra of infinitesimal gauge 
transformations. We have a natural homomorphism $\pi_* : \frak e 
\to \frak g$ with the kernel $\frak h$, i.e. $\frak e$ is an 
extension of $\frak g$ by means $\frak h$. 
 
Note that we have an additional structure of $C^\infty(M)$-module on 
$\frak g, \frak h, \frak e$, such that $[X,fY] = f[X,Y] + (\pi_*X)fY$, 
where $X,Y \in \frak e, \, f \in C^\infty(M)$. 
In particular, $\frak h$ is a Lie algebra over $C^\infty(M)$. 
The extension $$ 0 \rightarrow \frak h \rightarrow \frak e 
\rightarrow \frak g \rightarrow 0 $$ is also an extension of 
$C^\infty(M)$-modules.   
   
Assume now that the section $s : \frak g \to \frak e$   is a 
homomorphism of $C^\infty(M)$-modules. Then it can be considered as a 
connection in the principal bundle $\pi$, and the $\frak h$-valued 2-form 
$\rho$ as its curvature. In this sense we interpret the 
constructions from section \nmb!{1} as follows. See 
\cit!{7}, section~11 for more background information. The analogy 
with differential geometry has also been noticed by \cit!{8} and 
\cit!{9}.

\subhead\nmb.{4}. Geometric interpretation \endsubhead
Note that \thetag{\nmb!{2}.2} looks like the Maurer-Cartan formula 
for the \idx{\it curvature} on principal bundles of differential geometry
$$
\rh = ds + \tfrac12[s,s]_\wedge,
$$
where for an arbitrary vector space $V$ the usual Chevalley 
differential is given by 
$$\gather
d: L^p_{\text{skew}}(\g;V) 
     \to L^{p+1}_{\text{skew}}(\g;V)\\
d\ph(X_0,\dots,X_p) 
= \sum_{i<j}(-1)^{i+j}\ph([X_i,X_j],X_0,\dots,\widehat{X_i},\dots,
        \widehat{X_j},\ldots,X_p)
\endgather$$
and where for a vector space $W$ and a Lie algebra $\frak f$ the 
$\Bbb N$-graded (super) Lie bracket $[\quad,\quad]_\wedge$ 
on $L^*_{\text{skew}}(W,\frak f)$ is given by 
$$
[\ph,\ps]_\wedge(X_1,\dots,X_{p+q})
= \frac 1{p!\,q!} \sum_{\si} \text{sign}(\si)\,
        [\ph(X_{\si1},\dots,X_{\si p}),
     \ps(X_{\si(p+1)},\dotsc)]_{\frak f}.
$$
Similarly formula \thetag{\nmb!{2}.3} reads as 
$$
\ad_\rh = d\al + \tfrac12[\al,\al]_\wedge.
$$
Thus we view $s$ as a \idx{\it connection} in the sense of 
a \idx{\it horizontal lift} of vector fields on the base of a bundle, 
and $\al$ as an \idx{\it induced connection}. Namely, for every 
$\der(\h)$-module $V$ we put 
$$\gather
\al_\wedge : L^p_{\text{skew}}(\g;V) 
     \to L^{p+1}_{\text{skew}}(\g;V)\\
\al_\wedge \ph(X_0,\dots,X_p) 
= \sum_{i=0}^p(-1)^i\al_{X_i}(\ph(X_0,\dots,\widehat{X_i},\dots,X_p)).
\endgather$$
Then we have the \idx{\it covariant exterior differential} (on the 
sections of an associated vector bundle)
$$
\de_\al: L^p_{\text{skew}}(\g;V) 
\to L^{p+1}_{\text{skew}}(\g;V),\qquad \de_\al\ph = \al_\wedge \ph + d\ph, 
\tag{\nmb!{3}.1}
$$
for which formula \thetag{\nmb!{2}.4} looks like the \idx{\it Bianchi 
identity} $\de_\al\rh=0$. Moreover one finds quickly that another 
well known result from differential geometry holds, namely
$$
\de_\al\de_\al(\ph) = [\rh,\ph]_\wedge,\quad 
     \ph\in L^p_{\text{skew}}(\g;\h).\tag{\nmb!{3}.2}
$$
If we change the linear section $s$ to $s'=s+b$ for linear 
$b:\g\to\h$, then we get 
$$\align
\al'_X &= \al_X + \ad^\h_{b(X)} \tag{\nmb!{3}.3}\\
\rh'(X,Y) &= \rh(X,Y) + \al_Xb(Y) -\al_Yb(X) - b([X,Y]) + [bX,bY] 
\tag{\nmb!{3}.4}\\
&= \rh(X,Y) + (\de_\al b)(X,Y) + [bX,bY].\\
\rh'&= \rh + \de_\al b + \tfrac12 [b,b]_\wedge. 
\endalign$$

\proclaim{\nmb.{5}. Theorem}
Let $\h$ and $\g$ be Lie algebras. 

Then isomorphism classes of extensions of $\g$ over $\h$, i.e\. short 
exact sequences of Lie algebras 
$0\to \h\to \e\to \g\to 0$, modulo the equivalence described by the 
commutative diagram of Lie algebra homomorphisms
$$\CD
0 @>>> \h @>>> \e @>>> \g @>>> 0 \\
@.     @|  @V{\ph}VV   @|  @. \\
0 @>>> \h @>>> \e' @>>> \g @>>> 0, \\
\endCD$$
correspond bijectively
to equivalence classes of data of the following form:
$$\align
&\text{A linear mapping }
     \al:\g\to \der(\h),\tag{\nmb!{5}.1}\\
&\text{a skew-symmetric bilinear mapping }
     \rh:\g\x\g\to \h \tag{\nmb!{5}.2}\\ 
\endalign$$
such that
$$\align
&[\al_X,\al_Y]-\al_{[X,Y]} = \ad_{\rh(X,Y)}, \tag{\nmb!{5}.3}\\
&\sum_{\text{cyclic}} \Bigl(\al_X \rh(Y,Z) - \rh([X,Y],Z) \Bigr) 
     =0\quad\text{ equivalently, }\de_\al\rh=0.\tag{\nmb!{5}.4}
\endalign$$
On the vector space $\e:=\h\oplus \g$ a Lie algebra structure is 
given by 
$$
[H_1+X_1,H_2+X_2]_\e = [H_1,H_2]_\h + \al_{X_1}H_2 - \al_{X_2}H_1 + 
\rh(X_1,X_2) + [X_1,X_2]_\g,\tag{\nmb!{5}.5}
$$
the associated exact sequence is  
$$
0\to \h @>i_1>>\h\oplus\g= \e @>{\operatorname{pr}_2}>> \g\to 0.$$
Two data $(\al,\rh)$ and $(\al',\rh')$ are equivalent if there exists
a linear mapping $b\colon\g\to\h$ such that 
$$\align
\al'_X &= \al_X + \ad^\h_{b(X)}, \tag{\nmb!{5}.6}\\
\rh'(X,Y) &= \rh(X,Y) + \al_Xb(Y) -\al_Yb(X) - b([X,Y]) + [b(X),b(Y)] 
\tag{\nmb!{5}.7}\\
\rh'&= \rh + \de_\al b + \tfrac12 [b,b]_\wedge,\\
\endalign$$
the corresponding isomorphism being
$$
\e=\h\oplus\g \to \h\oplus\g=\e',\qquad H+X\mapsto H-b(X)+ X.
$$
Moreover, a datum $(\al,\rh)$ corresponds to a split extension (a 
semidirect product) if and only if $(\al,\rh)$ is equivalent to to a 
datum of the form $(\al',0)$ (then $\al'$ is a homomorphism). This is 
the case if and only if there exists a mapping $b:\g\to \h$ such that
$$
\rh = -\de_\al b - \tfrac12[b,b]_\wedge.\tag{\nmb!{5}.8} 
$$
\endproclaim
\demo{Proof}
Straigthforward computations. 
\qed\enddemo

\proclaim{\nmb.{6}. Corollary} \cit!{10}
Let $\g$ and $\h$ be Lie algebras such that $\h$ has no center.
Then isomorphism classes of extensions of $\g$ over $\h$ correspond 
bijectively to Lie homomorphisms 
$$
\bar\al:\g\to \operatorname{out}(\h)=\der(\h)/\ad(\h). 
$$
\endproclaim

\demo{Proof} 
If $(\al,\rh)$ is a data, then the map
$\bar\al: \frak g \to \der(\h)/\ad(\h)$
is a Lie algebra homomorphism by \thetag{\nmb!{5}.3}. 
Conversely, let $\bar\al$ be given. 
Choose a linear lift $\al:\g\to \der(\h)$ of 
$\bar\al$.
Since $\bar\al$ is a Lie algebra homomorphism and $\h$ has no center, 
there is a uniquely defined skew symmetric linear mapping
$\rh:\g\x\g\to \h$ such that
$[\al_X,\al_Y]-\al_{[X,Y]} = \ad_{\rh(X,Y)}$.
Condition \thetag{\nmb!{5}.4} is then automatically 
satisfied. For later use also, we record the simple proof:
$$\align
&\sum_{\text{cyclic} X,Y,Z}\Bigl[\al_X \rh(Y,Z) 
     - \rh([X,Y],Z),H \Bigr] \\
&=\sum_{\text{cyclic} X,Y,Z}\Bigl(\al_X [\rh(Y,Z),H] 
     - [\rh(Y,Z),\al_XH] - [\rh([X,Y],Z),H] \Bigr) \\
&=\sum_{\text{cyclic} X,Y,Z}\Bigl(\al_X [\al_Y,\al_Z] 
     - \al_X\al_{[Y,Z]} - [\al_Y,\al_Z]\al_X 
     + \al_{[Y,Z]}\al_X \\
&\qquad\qquad\qquad\qquad\qquad\qquad\qquad\qquad\qquad\qquad
     - [\al_{[X,Y]},\al_Z] + \al_{[[X,Y]Z]} \Bigr)H \\
&=\sum_{\text{cyclic} X,Y,Z}\Bigl([\al_X,[\al_Y,\al_Z]] 
     -[\al_X,\al_{[Y,Z]}] - [\al_{[X,Y]},\al_Z] 
     + \al_{[[X,Y]Z]} \Bigr)H =0.\\
\endalign$$
Thus $(\al,\rh)$ describes an extension by theorem \nmb!{5}. 
The rest is clear.
\qed\enddemo

\subhead\nmb.{7}. Remarks \endsubhead
If $\h$  has no center and 
$\bar\al:\g\to \operatorname{out}(\h)=\der(\h)/\ad(\h)$ is a given 
homomorphism, the  
extension corresponding to $\bar\al$ can be constructed in the 
following easy way: It is given by the pullback diagram
$$\CD
0 @>>>\h @>>>\der(\h) \x_{\operatorname{out}(\h)} \g @>{\pr_2}>>\g @>>> 0\\
@.     @|       @V{\pr_1}VV                            @V{\bar\al}VV @.\\
0 @>>> \h @>>> \der(\h) @>{\pi}>> \operatorname{out}(\h) @>>> 0\\
\endCD$$
where $\der(\h) \x_{\operatorname{out}(\h)} \g$ is the Lie subalgebra
$$
\der(\h) \x_{\operatorname{out}(\h)} \g 
     := \{(D,X)\in \der(\h)\x \g:\pi(D)=\bar\al(X)\} \subset 
     \der(\h)\x\g.
$$
We owe this remark to E\. Vinberg.

If $\h$ has no center and satisfies $\der(\h)=\h$, and if $\h$ is 
normal in a Lie algebra $\e$, then $\e\cong \h\oplus\e/\h$, since 
$\operatorname{Out}(\h)=0$.

\proclaim{\nmb.{8}. Theorem}
Let $\g$ and $\h$ be Lie algebras and let 
$$\bar\al:\g\to \operatorname{out}(\h)=\der(\h)/\ad(\h)$$
be a Lie algebra homomorphism. Then the following are equivalent:
\roster
\item For one (equivalently: any) linear lift 
       $\al:\g\to \der(\h)$ of $\bar\al$ choose 
       $\rh:\bigwedge^2\g\to \h$ satisfying 
       $([\al_X,\al_Y]-\al_{[X,Y]})=\ad_{\rh(X,Y)}$. Then the 
       $\de_{\bar\al}$-cohomology class of 
       $\la=\la(\al,\rh):=\de_{\al}\rh:\bigwedge^3\g \to Z(\h)$ in 
       $H^3(\g;Z(\h))$ vanishes. 
\item There exists an extension $0\to \h\to\e\to\g\to 0$ inducing the 
       homomorphism $\bar\al$. 
\endroster
If this is the case then all extensions $0\to \h\to\e\to\g\to 0$ 
inducing the homomorphism $\bar\al$ are parameterized by 
$H^2(\g,(Z(\h),\bar\al))$, the second Chevalley cohomology space of 
$\g$ with values in the center $Z(\h)$, considered as $\g$-module via 
$\bar\al$.
\endproclaim
\demo{Proof}
Using once more the computation in the proof of corollary \nmb!{6} we 
see that 
$\ad(\la(X,Y,Z))=\ad(\de_{\al}\rh(X,Y,Z))=0$ so that 
$\la(X,Y,Z)\in Z(\h)$.
The Lie algebra $\operatorname{out}(\h)=\der(\h)/\ad(\h)$ 
acts on the center $Z(\h)$, thus $Z(\h)$ is a $\g$-module via 
$\bar \al$, and $\de_{\bar\al}$ is the differential of the Chevalley 
cohomology. 
Using \thetag{\nmb!{3}.2} we see that
$$
\de_{\bar\al}\la = \de_{\al}\de_{\al}\rh = 
     [\rh,\rh]_\wedge = -(-1)^{2\cdot 2}[\rh,\rh]_\wedge =0,
$$
so that $[\la]\in H^3(\g;Z(\h))$. 

Let us check next that the cohomology class $[\la]$ does not depend 
on the choices we made. If we are given a pair $(\al,\rh)$ as above 
and we take another linear lift $\al':\g\to\der(\h)$ then 
$\al'_X=\al_X+\ad_{b(X)}$ for some linear $b:\g\to\h$. We consider 
$$
\rh':\bigwedge^2\g\to\h,\quad \rh'(X,Y)=\rh(X,Y)+(\de_\al b)(X,Y)+[b(X),b(Y)].
$$
Easy computations show that
$$\gather
[\al'_X,\al'_Y]-\al'_{[X,Y]} = \ad_{\rh'(X,Y)}\\
\la(\al,\rh) = \de_\al\rh = \de_{\al'}\rh' = \la(\al',\rh')
\endgather$$ 
so that even the cochain did not change. 
So let us consider for fixed $\al$ two linear mappings 
$$
\rh,\rh':\bigwedge^2\g\to \h,\quad 
[\al_X,\al_Y]-\al_{[X,Y]} = \ad_{\rh(X,Y)}= \ad_{\rh'(X,Y)}. 
$$
Then $\rh-\rh'=:\mu:\bigwedge^2\g\to Z(\h)$ and clearly
$\la(\al,\rh)-\la(\al,\rh')=\de_\al\rh-\de_\al\rh'=\de_{\bar\al}\mu$.

If there exists an extension inducing $\bar\al$ then for any lift 
$\al$ we may find $\rh$ as in \nmb!{5} such that $\la(\al,\rh)=0$.
On the other hand, given a pair $(\al,\rh)$ as in \therosteritem1 
such that $[\la(\al,\rh)]=0\in H^3(\g,(Z(\h),\bar\al))$, there 
exists $\mu:\bigwedge^2\g\to Z(\h)$ such that $\de_{\bar\al}\mu=\la$.
But then 
$$
\ad_{(\rh-\mu)(X,Y)}=\ad_{\rh(X,Y)},\quad \de_\al(\rh-\mu)=0,
$$ 
so that $(\al,\rh-\mu)$ satisfy the conditions of \nmb!{5} and thus 
define an extension which induces $\bar\al$.

Finally, suppose that \therosteritem1 is satisfied, and let us 
determine how many extensions there exist which induce $\bar\al$.
By \nmb!{5} we have to determine all equivalence classes of data 
$(\al,\rh)$ as in  
\nmb!{5}. We may fix the linear lift $\al$ and one mapping 
$\rh:\bigwedge^2\g\to \h$ which satisfies 
\thetag{\nmb!{5}.3} and \thetag{\nmb!{5}.4}, and we have to find all 
$\rh'$ with this property. But then 
$\rh-\rh'=\mu:\bigwedge^2\g\to Z(\h)$ and 
$$
\de_{\bar\al}\mu = \de_\al\rh-\de_\al\rh'=0-0=0
$$
so that $\mu$ is a 2-cocycle. Moreover we may still 
pass to equivalent data in the sense of \nmb!{5} using some 
$b:\g\to\h$ which does not change $\al$, i.e\. $b:\g\to Z(\h)$. The 
corresponding $\rh'$ is, by \thetag{\nmb!{5}.7},
$\rh'=\rh+\de_\al b + \tfrac12[b,b]_{\wedge }=\rh+\de_{\bar\al}b$.
Thus only the cohomology class of $\mu$ matters.  
\qed\enddemo

\proclaim{\nmb.{9}. Corollary}
Let $\g$ and $\h$ be Lie algebras such that $\h$ is abelian.
Then isomorphism classes of extensions of $\g$ over $\h$ correspond 
bijectively to the set of all pairs $(\al,[\rh])$, where
$\al:\g\to \g\frak l(\h)=\der(\h)$ is a homomorphism of Lie algebras 
and $[\rh]\in H^2(\g,\h)$ is a Chevalley cohomology class with 
coefficients in the $\g$-module $\h$. 
\endproclaim

\demo{Proof} This is obvious from theorem \nmb!{8}. 
\qed\enddemo

\subhead\nmb.{10}. An interpretation of the class $\la$ \endsubhead
Let $\h$ and $\g$ be Lie algebras and let a homomorphism 
$\bar\al:\g\to \der(\h)/\ad(\h)$ be given. We consider the extension 
$$
0\to \ad(\h) \to \der(\h) \to \der(\h)/\ad(\h) \to 0
$$
and the following diagram, where the bottom right hand square is a 
pullback (compare with remark \nmb!{7}):
$$\cgaps{0.5;0.5;0.5;1;1;0.5}\rgaps{0.5;1;1;0.4;1.2}\newCD
&&         & 0 @(0,-1) & 0 @(0,-1) & & \\
&&         & Z(\h) @()\a=@(1,0) @(0,-1) & Z(\h) @()\a-@(0,-1)  &    &  \\
&& 0 @(1,0) & \h @()\a-@(1,0) @(0,-1)    
     & \e @()\a-@(0,-1) @()\a-@(1,0) &  \g @(1,0) @()\a=@(0,-1) & 0  \\
&& 0 @(1,0) 
   & \ad(\h) @()\L{i}@(1,0) @()\a=\ds(4;0)\dtX(-4;0)@(-2,-2) @(0,-1) 
     & \e_0 @()\L{\be}\l{\text{\quad pull back}}@(-1,-2) @()\L{p}@(1,0) @(0,-1)
       & \g @(1,0) @()\dtX(-10;0)\l{\bar\al}@(-1,-2) & 0\\
&&         & 0 & 0 &  & \\
0 @(1,0) & \ad(\h) @(2,0) &  
     & \der(\h) @(1,0) 
             &  \der(\h)/\ad(\h) @(1,0) & 0 & \\
\endnewCD$$

The left hand vertical column describes $\h$ as a central extension of 
$\ad(\h)$ with abelian kernel $Z(\h)$ which is moreover killed 
under the action of $\g$ via $\bar\al$; it is given by a cohomology 
class $[\nu]\in H^2(\ad(\h);Z(\h))^\g$. 
In order to get an extension $\e$ of $\g$  
with kernel $\h$ as in the third row we have to check that the 
cohomology class $[\nu]$ is in the image of 
$i^*:H^2(\e_0;Z(\h))\to H^2(\ad(\h);Z(\h))^\g$.
It would be interesting to express this in terms of 
of the Hochschild-Serre 
exact sequence, see \cit!{6}.


\Refs

\widestnumber\key{44}

\ref
\key \cit0{1} 
\by Azc\'arraga, Jos\'e A.; Izquierdo, Jos\'e M.
\book Lie groups, Lie algebras, cohomology and some applications in 
physics
\bookinfo Cambridge Monographs on Mathematical Physics
\publ Cambridge University Press
\publaddr Cambridge, UK
\yr 1995
\endref

\ref
\key \cit0{2}
\by Eilenberg, S.; MacLane, S.
\paper Cohomology theory in abstract groups, II. Groups extensions 
with non-abelian kernel
\jour Ann. Math. (2)
\vol 48
\yr 1947
\pages 326--341
\endref

\ref
\key \cit0{3}
\by Giraud, Jean
\book Cohomologie non ab\'elienne
\bookinfo Grundlehren 179
\publ Springer-Verlag
\publaddr Berlin etc.
\yr 1971
\endref

\ref
\key \cit0{4}
\by Hochschild, G.
\paper Cohomology clases of finite type and finite dimensional 
kernels for Lie algebras 
\jour Am. J. Math.
\vol 76
\yr 1954
\pages 763-778
\endref

\ref
\key \cit0{5}
\by Hochschild, G. P.; Serre, J.-P.
\paper Cohomology of group extensions
\jour Trans. AMS 
\vol 74
\yr 1953
\pages 110--134
\endref

\ref
\key \cit0{6}
\by Hochschild, G. P.; Serre, J.-P.
\paper Cohomology of Lie algebras
\jour Ann. Math.
\vol 57
\yr 1953
\pages 591--603
\endref

\ref 
\key \cit0{7}
\by Kol\'a\v r, I.; Michor, Peter W.; Slov\'ak, J. 
\book Natural operations in differential geometry 
\publ Springer-Verlag
\publaddr Berlin Heidelberg New York
\yr 1993
\endref

\ref
\key \cit0{8}
\by Lecomte, Pierre
\paper Sur la suite exacte canonique asoci\'ee \`a un fibr\'e principal
\jour Bul. Soc. Math. France
\vol 13
\yr 1985
\pages 259--271
\endref

\ref
\key \cit0{9}
\by Lecomte, P.
\paper On some Sequence of graded Lie algebras asociated to manifolds
\jour Ann. Global Analysis Geom. 
\vol 12
\yr 1994
\pages 183--192
\endref

\ref
\key \cit0{10}
\by Lecomte, P.; Roger, C.
\paper Sur les d\'eformations des alg\`ebres de courants de type r\'eductif
\jour C. R. Acad. Sci. Paris, I
\vol 303
\yr 1986
\pages 807--810 
\endref

\ref
\key \cit0{11}
\by Mackenzie, K.
\book Lie groupoids and Lie algebroids in diferential geometry
\bookinfo London Mathematical Society Lecture Note Series, 124
\publ Cambridge University Press
\yr 1987
\endref

\ref
\key \cit0{12}
\by Mori, Mitsuya
\paper On the thre-dimensional cohomology group of Lie algebras
\jour J. Math. Soc. Japan
\vol 5
\yr 1953
\pages 171-183
\endref

\ref
\key \cit0{13}
\by Serre, J.-P.
\paper Cohomologie des groupes discrets
\jour Ann. of Math. Studies
\vol 70
\yr 1971
\pages 77--169
\finalinfo Princeton University Press
\endref

\ref
\key \cit0{14}
\by Serre, J.-P.
\paper Cohomologie des groupes discrets
\jour S\'eminaire Bourbaki
\vol 399
\yr 1970/71
\endref

\ref
\key \cit0{15}
\by Shukla, U.
\paper A cohomology for Lie algebras
\jour J. Math. Soc. Japan
\vol 18
\yr 1966
\pages 275-289
\endref

\endRefs
\enddocument